# A Class of Extended Ishikawa Iterative Processes in Banach Spaces for Nonexpansive Mappings

Congdian Cheng, Hongyan Guan
(College of Mathematics and Systems Science,
Shenyang Normal University, Shenyang, 110034, China)

**Abstract:** A class of extended Ishikowa Iterative processes is proposed and studied, which involves many kinds of Mann and Ishikawa iterative processes. The main conclusion of the present work extends and generalizes some recent results of this research line.

## 1  Introduction and Preliminaries

Throughout the present work, we always assume that $E$ is a real Banach space, $E^*$ is the duality space of $E$, $<\cdot,\cdot>$ is the dual pair between $E$ and $E^*$, and $J: E \to 2^{E^*}$ is the normalized duality mapping defined by

$$J(x) = \{f \in E^* : <x, f> = \|x\| \cdot \|f\|, \|f\| = \|x\|\}, x \in E. \qquad (1.1)$$

We also assume that $D$ is a nonempty closed convex subset of $E$, $T : D \to D$ is a mapping, and $F(T)$ denotes the set of all the fixed points of $T$. In addition, we use $\Pi_D$ representing all the contractions on $D$, i.e. $\Pi_D = \{f \mid f : D \to D,$ and there exists $\alpha \in (0,1)$ such that $\|f(x) - f(y)\| \le \alpha \|x - y\|$ for all the $x, y \in D\}$.

**Definition 1.** Let $U = \{x \in E : \|x\| = 1\}$. $E$ is said to be uniformly smooth if the limit $\lim_{t \to 0} \frac{\|x + ty\| - \|x\|}{t}$ is attained uniformly for all the $(x, y) \in U \times U$.

The following result is well known (see Goebel and Reich [1]).

**Proposition 1.** Let $E$ be uniformly smooth. Then the duality mapping $J$ defined by (1.1) is single valued, and it is uniformly continuous on the bounded subset of $E$ with the norm topologies of $E$ and $E^*$.

Recall that the sequences

$$x_{n+1} = \alpha_n u + (1 - \alpha_n) T x_n \ (or\ \alpha_n x_n + (1 - \alpha_n) T x_n), \qquad (1.2)$$

$$\begin{cases} x_{n+1} = \alpha_n u + (1 - \alpha_n) T y_n \ (or\ \alpha_n x_n + (1 - \alpha_n) T y_n) \\ y_n = \beta_n x_n + (1 - \beta_n) T x_n, \end{cases} \qquad (1.3)$$

$$x_{n+1} = \alpha_n u + (1-\alpha_n - \gamma_n)Tx_n + \gamma_n u_n, \quad (1.4)$$

$$(or\ x_{n+1} = \alpha_n x_n + (1-\alpha_n - \gamma_n)Tx_n + \gamma_n u_n)$$

$$\begin{cases} x_{n+1} = \alpha_n u + (1-\alpha_n - \gamma_n)Ty_n + \gamma_n u_n \\ (or\ x_{n+1} = \alpha_n x_n + (1-\alpha_n - \gamma_n)Ty_n + \gamma_n u_n) \\ y_n = \beta_n x_n + (1-\beta_n - \delta_n)Tx_n + \delta_n v_n, \end{cases} \quad (1.5)$$

are respectively called Mann iterative process, Ishikawa iterative process, modified Mann iterative process with error and modified Ishikawa iterative process with errors of $T$, where $x_0, u \in D$ and $n \geq 0$. The problem for these iterative sequences converging to the fixed point of $T$ was studied by lots of authors, e.g., Halpern [2], Reich [3], Zhang and Tian [4], Chidume [5], Liu [6], Liu Q.H. and Liu Y. [7], Zhao and Zhang [8]. In particular, Xu[9] generalized (1.2) to the iterative process

$$x_{n+1} = \alpha_n f(x_n) + (1-\alpha_n)Tx_n. \quad (1.6)$$

Under a certain conditions, he proved that $\{x_n\}$ converges strongly to a fixed point of $T$ and other related results. In 2007, Chang [10] extended and improved the work of Xu [9].

Motivated and inspired by the contributions above, the present work addresses the following iterative process.

$$\begin{cases} x_0 \in D, \quad n \geq 0 \\ x_{n+1} = \alpha_n f_n(x_n) + (1-\alpha_n)Ty_n \\ y_n = \beta_n g_n(x_n) + (1-\beta_n)Tx_n, \end{cases} \quad (1.7)$$

where $f_n$ and $g_n$ are all the mappings of $D$ into its self. To some extent, (1.7) generalizes and unifies all the iterative sequences mentioned above. For example, if we instead $\alpha_n, \beta_n, f_n$ and $g_n$ with $(\alpha_n + \gamma_n), (\beta_n + \delta_n), (\frac{\alpha_n u}{\alpha_n + \gamma_n} + \frac{\gamma_n u_n}{\alpha_n + \gamma_n})$ and $(\frac{\beta_n x}{\beta_n + \delta_n} + \frac{\delta_n v_n}{\beta_n + \delta_n})$ in (1.7) respectively, then it reduces to (1.5). We still concern with the problem whether the sequence (1.7) converges to the fixed point of $T$. To prove our main results, we also need the following Lemmas.

**Lemma 1.1.**[9] Let $X$ be a uniformly smooth Banach space, $C$ be a closed convex subset of $X$, $T: C \to C$ be a nonexpansive with $F(T) \neq \phi$, and $f \in \Pi_C$. Then $\{x_t\}$ defined by $x_t = tf(x_t) + (1-t)Tx_t$ converges strongly to a point in $F(T)$. If we define $Q: \Pi_C \to F(T)$ by

$$Q(f) := \lim_{t \to 0} x_t, \quad f \in \Pi_C, \tag{1.8}$$

then $Q(f)$ solves the variational inequality

$$< (I-f)Q(f), J(Q(f)-p) > \le 0, \quad f \in \Pi_C, p \in F(T).$$

In particular, if $f = u \in C$ is a constant, then (1.8) is reduced to the sunny nonexpansive retraction of Reich from $C$ onto $F(T)$,

$$< Q(u)-u, J(Q(u)-p) > \le 0, \quad u \in C, p \in F(T).$$

**Lemma 1.2.**[11] Let $X$ be a real Banach space and $J_p : X \to 2^{X^*}, 1 < p < \infty$, be a duality mapping. Then for any given $x, y \in X$, we have

$$\| x+y \|^p \le \| x \|^p + p \cdot < y, j_p >, \forall j_p \in J_p(x+y).$$

**Lemma 1.3.**[6] Let $\{a_n\}, \{b_n\}$ and $\{c_n\}$ be three nonnegative real sequences satisfying

$$a_{n+1} \le (1-t_n)a_n + b_n + c_n$$

with $\{t_n\} \subset [0,1], \sum_{n=0}^{\infty} t_n = \infty, b_n = o(t_n)$, and $\sum_{n=0}^{\infty} c_n < \infty$. Then $\lim_{n \to \infty} a_n = 0$.

## 2  Main results

In this section, we address the strong convergence of the iterative sequence (1.7).

**Lemma 2.1.** Let $f, f_n \in \Pi_D, t_n \in (0,1)$, let $T$ be a nonexpansive mapping, and let $z_n$ be the unique solution of the equation $z = t_n f_n(z) + (1-t_n)Tz$ for all $n \ge 0$. Then $z \to Q(f)$ (defined by (1.8) ) as $t_n \to 0$ (strongly) if $\{f_n(x)\}$ converges uniformly to $f(x)$ on $D$.

**Proof.** Let $z_n'$ be the unique solution of equation $z = t_n f(z) + (1-t_n)Tz$. Then

$$z_n - z_n' = t_n[f_n(z_n) - f(z_n')] + (1-t_n)[Tz_n - Tz_n'].$$

This leads to

$$\| z_n - z_n' \| \le t_n (\| f_n(z_n) - f(z_n) \| + \| f(z_n) - f(z_n') \|) + (1-t_n) \| z_n - z_n' \|$$

$$\le t_n \| f_n(z_n) - f(z_n) \| + \alpha t_n \| z_n - z_n' \| + (1-t_n) \| z_n - z_n' \|.$$

Thus $t_n(1-\alpha)\|z_n - z_n'\| \leq t_n\|f_n(z_n) - f(z_n)\|$. That is,

$$(1-\alpha)\|z_n - z_n'\| \leq \|f_n(z_n) - f(z_n)\|.$$

Since $f_n(x)$ converges uniformly to $f(x)$, $\|z_n - z_n'\| \to 0$. On the other hand, $z_n' \to Q(f)$ as $t_n \to 0$ from Lemma 1.1. Hence, $z_n \to Q(f)$ as $t_n \to 0$. This completes the proof.

**Theorem 2.1.** For (1.7), let $T$ be a nonexpansive mapping with $F(T) \neq \phi$, and $\{f_n\}$ converges uniformly to $f \in \Pi_D$ and uniform contractions, that is, there exists $L \in (0,1)$ such that $\|f_n(x) - f_n(y)\| \leq L\|x - y\|$ for all $x, y \in D$, and for all the indices $n$. Moreover, suppose $\beta_n \|g_n(x) - x\| \leq \delta_n(\|x\| + M)$, where $M$ is a positive constant, and $\sum_{n=0}^{\infty} \delta_n < \infty$; $\alpha_n \to 0$, $\sum_{n=0}^{\infty} \alpha_n = \infty$, $\alpha_n \in [0,1]$; $\beta_n \in (0,1]$. In addition, assume also $E$ is uniformly smooth, $\{x_n\}$ is bounded and $\|Tx_n - x_n\| \to 0$. Then $\{x_n\}$ converges strongly to $Q(f)$, where $Q: \Pi_D \to F(T)$ is the defined by (1.8).

**Proof.** Since $\|Tx_n - x_n\| \to 0$, we can choose $\{t_n\}$ such that $\|Tx_n - x_n\| = o(t_n)$. Let $z_n$ be the unique solution of the fixed point of equation $z = t_n f_n(z) + (1-t_n)Tz$. Then $\{z_n\}$ converges strongly to $Q(f)$ by Lemma 2.1. Let $z = Q(f)$. Then we have

$$\|x_{n+1} - z\|^2 = \|(1-\alpha_n)(Ty_n - z) + \alpha_n(f_n(x_n) - z)\|^2$$

$$\leq (1-\alpha_n)^2 \|Ty_n - z\|^2 + 2\alpha_n <f_n(x_n) - z, j(x_{n+1} - z)> \quad \text{(by Lemma 1.2)}$$

$$\leq (1-\alpha_n)^2 \|y_n - z\|^2 + 2\alpha_n(<f_n(x_n) - f_n(z), j(x_{n+1} - z)>$$

$$+ <f_n(z) - z, j(x_{n+1} - z)>)$$

$$\leq (1-\alpha_n)^2 \|y_n - z\|^2 + 2\alpha_n L \|x_n - z\| \cdot \|x_{n+1} - z\|$$

$$+ 2\alpha_n <f_n(z) - z, j(x_{n+1} - z)>$$

$$\leq (1-\alpha_n)^2 \|y_n - z\|^2 + \alpha_n L(\|x_n - z\|^2 + \|x_{n+1} - z\|^2)$$

$$+ 2\alpha_n < f_n(z) - z, j(x_{n+1} - z) >; \tag{2.1}$$

$$\| y_n - z \| = \| \beta_n(g_n(x_n) - x_n) + \beta_n(x_n - z) + (1 - \beta_n)(Tx_n - Tz) \|$$

$$\leq \delta_n(\| x_n \| + M) + \| x_n - z \|,$$

$$\| y_n - z \|^2 \leq [\delta_n(\|x_n\| + M)]^2 + 2[\delta_n(\|x_n\| + M)]\|x_n - z\| + \|x_n - z\|^2$$
$$= \delta_n[\delta_n(\|x_n\| + M)^2 + 2(\|x_n\| + M)\|x_n - z\|] + \|x_n - z\|^2 \tag{2.2}$$
$$\leq \delta_n M_1 + \| x_n - z \|^2.$$

On the other hand, we have

$$\| z_n - x_n \|^2 = \|(1 - t_n)(Tz_n - x_n) + t_n(f_n(z_n) - x_n)\|^2$$

$$\leq (1 - t_n)^2 \| Tz_n - x_n \|^2 + 2t_n(< f_n(z_n) - z_n, j(z_n - x_n) >$$

$$+ < z_n - x_n, j(z_n - x_n) >)$$

$$= (1 - t_n)^2 \| Tz_n - x_n \|^2 + 2t_n(< f_n(z_n) - z_n, j(z_n - x_n) - j(z - x_{n+1}) >$$

$$+ < f_n(z_n) - z_n, j(z - x_{n+1}) > + \| z_n - x_n \|^2), \tag{2.3}$$

$$\| Tz_n - x_n \|^2 = \| Tz_n - Tx_n + Tx_n - x_n \|^2 \leq (\| Tz_n - Tx \| + \| Tx_n - x_n \|)^2$$

$$= \| Tz_n - Tx_n \|^2 + [(2\| Tz_n - Tx_n \| + \| Tx_n - x_n \|)] \cdot \| Tx_n - x_n \|$$

$$\leq \| z_n - x_n \|^2 + M_2 o(t_n), \tag{2.4}$$

(Note: $z_n \to z$ and $\{x_n\}$ is bounded.) and

$$< f_n(z_n) - z_n, j(z_n - x_n) - j(z - x_{n+1}) >$$

$$\leq \| f_n(z_n) - z_n \| \cdot \| j(z_n - x_n) - j(z - x_{n+1}) \|$$

$$\leq M_3(\| j(z_n - x_n) - j(z - x_n) \| + \| j(z - x_n) - j(z - x_{n+1}) \|), \tag{2.5}$$

(Note: $\| f_n(z_n) - z_n \| = \| f_n(z_n) - f(z_n) + f(z_n) - f(z') + f(z') - z_n \|$

$\leq \| f_n(z_n) - f(z_n) \| + (1 + \alpha) \| z' - z_n \| \leq M_3$, where $z'$ is a fixed point of $f$.)

$$< f_n(z_n) - z_n, j(z - x_{n+1}) >$$

$$= < [f_n(z_n) - f_n(z)] + [f_n(z) - z] + (z - z_n), j(z - x_{n+1}) >$$

$$\leq (1 + L) \| z - z_n \| \cdot \| z - x_{n+1} \| + < f_n(z) - z, j(z - x_{n+1}) >$$

$$= M_4 \| z - z_n \| + < f_n(z) - z, j(z - x_{n+1}) > . \tag{2.6}$$

Substitute in (2.3) the (2.4), (2.5) and (2.6), we obtain

$$< z - f_n(z), j(z - x_{n+1}) > \leq \frac{t_n}{2} \| z_n - x_n \|^2 + \frac{M_2 o(t_n)}{2t_n} + M_3 (\| j(z_n - x_n)$$

$$- j(z - x_n) \| + \| j(z - x_n) - j(z - x_{n+1}) \|) + M_4 \| z - z_n \| = b_n \to 0 \ (n \to \infty). \tag{2.7}$$

(Note: In terms of Proposition 1, $j$ is uniformly continuous on bounded subset.) Combining (2.1), (2.2) and (2.7), we also obtain

$$(1 - L\alpha_n) \| x_{n+1} - z \|^2 \leq [1 - (2 - L)\alpha_n + \alpha_n^2] \cdot \| x_n - z \|^2 + 2\alpha_n b_n + M_1 \delta_n.$$

This further leads to

$$\| x_{n+1} - z \|^2 \leq [1 - \frac{2(1-L)\alpha_n + \alpha_n^2}{1 - L\alpha_n}] \cdot \| x_n - z \|^2 + \frac{2\alpha_n b_n}{1 - L\alpha_n} + \frac{M_1 \delta_n}{1 - L\alpha_n}.$$

By Lemma 1.2, $x_{n+1} \to z$. This completes the proof.

## 3 Special cases

When $\beta_n = 1$, $g_n(x) = x$ and $f_n(x) = f(x)$, (1.7) reduces to (1.6), and that $\{f_n(x)\}$ converges uniformly to $f$ holds obviously. By Theorem 2.1, we can immediately obtain the following conclusion, which is the major conclusion of [10, Theorem 1].

**Theorem 3.1.** For (1.6), let $T$ be a nonexpansive mapping with $F(T) \neq \phi$; $f \in \Pi_D, x_0 \in D$; $\alpha_n \to 0$, and $\sum_{n=0}^{\infty} \alpha_n = \infty$. Assume $\| Tx_n - x_n \| \to 0$. Then $\{x_n\}$ converges strongly to $Q(f)$ if $E$ be uniformly smooth. ( Note: In the case, $\{x_n\}$ is bounded from the Proof of [10, Theorem 1]. )

For (1.5), by putting $\alpha'_n = \alpha_n + \gamma_n$, $\beta'_n = \beta_n + \delta_n$, $f_n(x) = \frac{\alpha_n u}{\alpha_n + \gamma_n} + \frac{\gamma_n u_n}{\alpha_n + \gamma_n}$ and $g_n(x) = \frac{\beta_n x}{\beta_n + \delta_n} + \frac{\delta_n v_n}{\beta_n + \delta_n}$, it becomes (1.7) with $\alpha'_n$ and $\beta'_n$. We can easily derive the next Theorem 3.2 from Theorem 2.1, which roughly resembles the conclusion of [8, Theorem 2.3] in the uniformly smooth space.

**Theorem 3.2.** For (1.5), let $T$ be a nonexpansive mapping with $F(T) \neq \phi$; $0 \leq \alpha_n + \gamma_n \leq 1$, $0 < \beta_n + \delta_n \leq 1$; $\alpha_n \to 0$, and $\sum_{n=0}^{\infty} \alpha_n = \infty$; $\sum_{n=0}^{\infty} \gamma_n < \infty$, $\sum_{n=0}^{\infty} \delta_n < \infty$;

$\{u_n\}$ and $\{v_n\}$ are bounded in $D$; and $\dfrac{\gamma_n}{\alpha_n + \gamma_n} \to 0$. Assume also $E$ be uniformly smooth, $\{x_n\}$ is bounded and $\|Tx_n - x_n\| \to 0$. Then $\{x_n\}$ converges strongly to $Q(f)$.

**Proof.** Let $f_n(x) = \dfrac{\alpha_n u}{\alpha_n + \gamma_n} + \dfrac{\gamma_n u_n}{\alpha_n + \gamma_n}$, $g_n(x) = \dfrac{\beta_n x}{\beta_n + \delta_n} + \dfrac{\delta_n v_n}{\beta_n + \delta_n}$, $\alpha'_n = \alpha_n + \gamma_n$ and $\beta'_n = \beta_n + \delta_n$. It is obvious that $\{f_n(x)\}$ is uniform contractions, $\alpha'_n \in [0,1]$, $\beta'_n \in (0,1]$.

Since $\dfrac{\gamma_n}{\alpha_n + \gamma_n} \to 0$ and $\{u_n\}$ is bounded, $f_n(x) = u - \dfrac{\gamma_n}{\alpha_n + \gamma_n}(u - u_n) \to u = f(x) \in \Pi_D$. Since $\{v_n\}$ is bounded, $\beta'_n \|g_n(x) - x\| = \beta'_n \|\dfrac{\beta_n x}{\beta_n + \delta_n} + \dfrac{\delta_n v_n}{\beta_n + \delta_n} - x\| = \delta_n \|-x + v_n\| \leq \delta_n(\|x\| + M)$. Moreover, since $\alpha_n \to 0$, $\sum_{n=0}^{\infty} \alpha_n = \infty$ and $\sum_{n=0}^{\infty} \gamma_n < \infty$, we have $\alpha'_n \to 0$ and $\sum_{n=0}^{\infty} \alpha'_n = \infty$. Thus, Theorem 3.2 holds from Theorem 2.1.

In addition, Yao [12] also studied the sequence

$$\begin{cases} x_{n+1} = \alpha_n u + \beta_n x_n + \gamma_n T x_n \\ \alpha_n + \beta_n + \gamma_n = 1 \end{cases}, \quad (3.1)$$

which can be transformed as

$$x_{n+1} = (\alpha_n + \beta_n)(\dfrac{\beta_n x_n}{\alpha_n + \beta_n} + \dfrac{\alpha_n u}{\alpha_n + \beta_n}) + \gamma_n T x_n$$

$$= (\alpha_n + \beta_n)[x_n - \dfrac{\alpha_n}{\alpha_n + \beta_n}(x_n - u)] + \gamma_n T x_n.$$

Thus, we can easily know the following conclusion holds from Theorem 2.1, which can be taken as a complementary result of [12, Theorem 3.1].

**Theorem 3.3.** For (3.1), let $T$ be a nonexpansive mapping with $F(T) \neq \phi$, $x_0, u \in D$; $\alpha_n, \beta_n, \gamma_n \geq 0$; $\alpha_n \to 0$, $\beta_n \to 0$ and $\sum_{n=0}^{\infty} \alpha_n = \infty$. Then $\{x_n\}$ converges strongly to a point of $F(T)$.

References
[1] Goebel K., Reich S., Uniform convexity, nonexpansive mappings and hyperbolic geometry,